\documentclass[a4paper,11pt]{amsart}

\usepackage{graphicx}
\usepackage{mathptmx}
\usepackage{amsmath}
\usepackage{amssymb}
\usepackage{enumitem}
\usepackage{xcolor}

\newmuskip\pFqmuskip

\newcommand*\pFq[6][8]{%
  \begingroup % only local assignments
  \pFqmuskip=#1mu\relax
  \mathcode`=\string"8000
  \begingroup\lccode`\~=`\,
  \lowercase{\endgroup\let~}\pFqcomma
  F^{#2}_{#3}{\left(\genfrac..{0pt}{}{#4}{#5}\bigg|#6\right)}%
  \endgroup
}
\newcommand{\pFqcomma}{\mskip\pFqmuskip}

\newtheorem{theorem}{Theorem}[section]

\begin{document}

\title[]{Probabilistic multi-Stirling numbers of the second kind and probabilistic multi-Lah numbers}

\author{Taekyun  Kim*}
\address{Department of Mathematics, Kwangwoon University, Seoul 139-701, Republic of Korea}
\email{tkkim@kw.ac.kr}
\author{Dae San  Kim}
\address{Department of Mathematics, Sogang University, Seoul 121-742, Republic of Korea}
\email{dskim@sogang.ac.kr}

\thanks { * is corresponding author.}
\begin{abstract}
Assume that the moment generating function of the random variable $Y$ exists in a neighborhood of the origin. We introduce the probabilistic multi-Stirling numbers of the second kind associated with $Y$ and the probabilistic multi-Lah numbers associated with $Y$, both of indices $(k_{1},k_{2},\dots,k_{r})$, by means of the multiple logarithm. Those numbers are respectively probabilistic extensions of the multi-Stirling numbers of the second kind and the multi-Lah numbers which, for $(k_{1},k_{2},\dots,k_{r})=(1,1,\dots,1)$, boil down respectively to the Stirling numbers of the second and the unsigned Lah numbers. The aim of this paper is to study some properties, related identities, recurrence relations and explicit expressions of those probabilistic extension numbers in connection with several other special numbers.
\end{abstract}

\maketitle

\markboth{\centerline{\scriptsize Probabilistic multi-Stirling numbers of the second kind and probabilistic multi-Lah numbers}}
{\centerline{\scriptsize T. Kim and D. S. Kim}}

\section{Introduction}
Assume that $Y$ is a random variable whose moment generating function exists in a neighborhood of the origin (see \eqref{6-1}). The aim of this paper is to study the probabilistic multi-Stirling numbers of the second kind ${n \brace k_{1},k_{2},\dots k_{r}}_{Y}$ associated with $Y$ (see \eqref{18}) and the probabilistic multi-Lah numbers $L_{Y}^{(k_{1},k_{2},\dots,k_{r})}(n,r)$ associated with $Y$ (see \eqref{26}). We investigate some properties, related identities, recurrence relations and explicit expressions for those numbers in connection with the multi-Bernoulli numbers  $B_{n}^{(k_{1},k_{2},\dots,k_{r})}$ (see \eqref{4}), the probabilistic Fubini polynomials of order $r$, $F_{n}^{(r,Y)}(y)$, associated with $Y$ (see \eqref{31}), the probabilistic Stirling numbers of the second kind ${n \brace r}_{Y}$ associated with $Y$ (see \eqref{14}), the unsigned multi-Stirling numbers of the first kind ${n \brack k_{1},k_{2},\dots,k_{r}}$ (see \eqref{8}), the unsigned Stirling numbers of first kind ${n \brack r}$ (see \eqref{5}) and the Stirling numbers of the second kind ${n \brace r}$ (see \eqref{6}). \par
We remark that ${n \brace k_{1},k_{2},\dots k_{r}}_{Y}$ is a probabilistic extension of the multi-Stirling number of the second kind ${n \brace k_{1},k_{2},\dots k_{r}}$ (see \eqref{10}), which reduces to the Stirling number of the second kind ${n \brace r}$, for $(k_{1},k_{2},\dots k_{r})=(1,1,\dots,1)$. We also note that $L_{Y}^{(k_{1},k_{2},\dots,k_{r})}(n,r)$ is a probabilistic extension of the multi-Lah number $L^{(k_{1},\dots,k_{r})}(n,r)$ (see \eqref{2}), which boils down to the unsigned Lah number $L(n,r)$ (see \eqref{3-1}), for $(k_{1},k_{2},\dots k_{r})=(1,1,\dots,1)$. In addition, ${n \brack k_{1},k_{2},\cdots,k_{r}}$ and $B_{n}^{(k_1,k_2,\dots,k_{r})}$ generalize respectively ${n \brack r}$ and $\frac{(-1)^{n}}{r!}B_{n}^{(r)}$, for $(k_{1},k_{2},\dots k_{r})=(1,1,\dots,1)$, where $B_{n}^{(r)}$ are the Bernoulli numbers of order $r$. The common feature of the six numbers ${n \brace k_{1},k_{2},\dots k_{r}}_{Y}$, $L_{Y}^{(k_{1},k_{2},\dots,k_{r})}(n,r)$, ${n \brace k_{1},k_{2},\dots k_{r}}$, $L^{(k_{1},k_{2},\dots,k_{r})}(n,r)$, ${n \brack k_{1},k_{2},\cdots,k_{r}}$ and $B_{n}^{(k_1,k_2,\dots,k_{r})}$ is that they are all defined with the help of the multiple logarithm $\mathrm{Li}_{k_{1},k_{2},\dots,k_{r}}(t)$ (see \eqref{1}), which reduces to the polylogarithm $\mathrm{Li}_{k_{1}}(t)$, for $r=1$ and to $\frac{1}{r!}\big(-\log(1-t)\big)^{r}$, for $(k_{1},k_{2},\dots k_{r})=(1,1,\dots,1)$. \par

The outline of this paper is as follows. First, we recall the multiple logarithms, the multi-Lah numbers and the multi-Bernoulli numbers. We remind the reader of the unsigned Stirling numbers of the first kind and the Stirling numbers of the second kind. We recall the probabilistic Lah numbers $L_{Y}(n,k)$ associated with $Y$. We remind the reader of the unsigned multi-Stirling numbers of the first kind and the multi-Stirling numbers of the second kind. Finally, we recall the probabilistic Stirling numbers of the second. Section 2 contains the main results of this paper. We derive first that $\mathrm{Li}_{\underbrace{1,1,\dots,1}_{r-\mathrm{times}}}(t)=\frac{1}{r!}\big(-\log(1-t)\big)^{r}$. Then we define the probabilistic multi-Stirling numbers of the second kind ${n \brace k_{1},k_{2},\dots k_{r}}_{Y}$ associated with $Y$, and show ${n \brace \underbrace{1,1,\dots,1}_{r-\mathrm{times}}}\frac{}{}_{Y}={n \brace r}_{Y}$. In Theorem 2.1, we express ${n+1 \brace k_{1},k_{2},\dots,k_{r-1},1}$ as a finite sum involving ${m \brace k_{1},k_{2},\dots, k_{r-1}},\,\, (m=r-1,\dots,n)$, which generalizes a well-known identity for the Stirling numbers of the second kind to that for the probabilistic Stirling numbers of the second kind associated with $Y$. In Theorem 2.2, we express $\sum_{m=0}^{n}B_{m}^{(k_{1},k_{2},\dots,k_{r})}{n \brace m}_{Y}$ as a finite sum involving the product ${i \brace j}_{Y}{n-i \brace k_{1},k_{2},\dots,k_{r}}_{Y}$. In Theorem 2.3, we represent ${n \brace k_{1},k_{2},\dots,k_{r}}_{Y}$ as a finite sum involving the triple product ${l \brace m}{n \brace l}_{Y}{m \brack k_{1},k_{2},\dots,k_{r}}$. We define the probabilistic multi-Lah numbers $L_{Y}^{(k_{1},k_{2},\dots,k_{r})}(n,r)$ associated with $Y$. In Theorem 2.4, we show that $L_{Y}^{(k_{1},k_{2},\dots,k_{r})}(n,r)= \sum_{k=r}^{n}{n \brace k_{1},k_{2},\dots,k_{r}}_{Y}{n \brack k}$. In Theorem 2.5, we express ${n \brace k_{1},k_{2},\dots,k_{r}}_{Y}$ as a finite sum involving the triple product ${l \brace r}{n \brace m+l}_{Y}B_{m}^{(k_{1},k_{2},\dots,k_{r})}$. In Theorem 2.6, we show the identity $\sum_{k=r}^{n}{n \brace k}_{Y}L^{(k_{1},k_{2},\dots,k_{r})}(k,r)=\sum_{k=r}^{n}\binom{n}{k}{k \brace k_{1},k_{2},\dots,k_{r}}_{Y}F_{n-k}^{(r,Y)}(1)$, where $F_{n}^{(r,Y)}(y)$ are the probabilistic Fubini polynomials of order $r$. For the rest of this section, we recall the facts that are needed throughout this paper.

\vspace{0.1in}

For integers $k_{i}\ (1\le i\le r)$, the multiple logarithm of index $(k_{1},k_{2},\dots,k_{r})$ is given by
\begin{equation}
\mathrm{Li}_{k_{1},k_{2},\dots,k_{r}}(t)=\sum_{0<m_{1}<m_{2}<\cdots<m_{r}}\frac{t^{m_{r}}}{m_{1}^{k_{1}}m_{2}^{k_{2}}\cdots m_{r}^{k_{r}}},\quad (|t|<1).\quad (\mathrm{see}\ [19,24]). \label{1}
\end{equation}
In particular, for $r=1$, we have
\begin{displaymath}
\mathrm{Li}_{k_{1}}(t)=\sum_{m=1}^{\infty}\frac{t^{m}}{m^{k_{1}}},\quad (\mathrm{see}\ [19,24,28]),
\end{displaymath}
where $\mathrm{Li}_{k_{1}}(t)$ is called the polylogarithm. \par
Recently, the multi-Lah number of index $(k_{1},k_{2},\dots,k_{r})$ is defined by
\begin{equation}
\frac{\mathrm{Li}_{k_{1},k_{2},\dots,k_{r}}(1-e^{-t})}{(1-t)^{r}}=\sum_{n=r}^{\infty}L^{(k_{1},\dots,k_{r})}(n,r)\frac{t^{n}}{n!},\quad (\mathrm{see}\ [19,24]). \label{2}
\end{equation}
Note that
\begin{equation}
\sum_{n=r}^{\infty}L^{\overbrace{(1,1,\dots,1)}^{r-times}}(n,r)\frac{t^{n}}{n!}=\frac{\mathrm{Li}_{1,1,\dots,1}(1-e^{-t})}{(1-t)^{r}}=\frac{1}{r!}\bigg(\frac{t}{1-t}\bigg)^{r}=\sum_{n=r}^{\infty}L(n,r)\frac{t^{n}}{n!}.\label{3}
\end{equation}
Here $L(n,r)$ are the unsigned Lah numbers given by
\begin{equation}
\langle x\rangle_{n}=\sum_{k=0}^{n}L(n,k)(x)_{k},\quad (n\ge 0),\quad (\mathrm{see}\ [12,17,18,20-25,27]), \label{3-1}
\end{equation}
where $\langle x\rangle_{0}=1,\ \langle x\rangle_{n}=x(x+1)\cdots (x+n-1),\quad(n \ge 1)$, and $(-1)^{n}\langle -x\rangle_{n}=(x)_{n}$. \par
In [19], Kim-Kim introduced the multi-Bernoulli number of index $(k_{1},k_{2},\dots,k_{r})$ defined by
\begin{equation}
\frac{1}{(1-e^{-t})^{r}} \mathrm{Li}_{k_{1},k_{2},\dots,k_{r}}(1-e^{-t})=\sum_{n=0}^{\infty}B_{n}^{(k_{1},k_{2},\dots,k_{r})}\frac{t^{n}}{n!}. \label{4}
\end{equation}
Note that $B_{n}^{\overbrace{(1,1,\dots,1)}^{r-times}}=\frac{(-1)^{n}}{r!}B_{n}^{(r)}$, where $B_{n}^{(r)}$ are the Bernoulli numbers of order $r$, given by
\begin{equation*}
\bigg(\frac{t}{e^{t}-1}\bigg)^{r}=\sum_{n=0}^{\infty}B_{n}^{(r)}\frac{t^{n}}{n!},\quad (\mathrm{see}\ [1-32]).
\end{equation*}
The unsigned Stirling numbers of the first kind are defined by
\begin{equation}
\langle x\rangle_{n}=\sum_{k=0}^{n}{n \brack k}x^{k},\quad (n,k\ge 0),\quad (\mathrm{see}\ [17,18,20]).\label{5}
\end{equation}
The Stirling numbers of the second kind are defined by
\begin{equation}
x^{n}=\sum_{k=0}^{n}{n \brace k}(x)_{k},\quad (n\ge 0),\quad (\mathrm{see}\ [8,12,29]). \label{6}
\end{equation}
Throughout this paper, we assume that $Y$ is random variable such that the moment generating function of $Y$,
\begin{equation}
E\big[e^{tY}\big]=\sum_{n=0}^{\infty}\frac{t^{n}}{n!}E\big[Y^{n}\big],\quad (|t|<r),\quad\textrm{exists for some $r>0$.} \label{6-1}
\end{equation}
Let $(Y_{j})_{j\ge 1}$ be a sequence of mutually independent copies of random variable $Y$, and let
\begin{displaymath}
S_{k}=Y_{1}+Y_{2}+\cdots+Y_{k},\ (k\ge 1),\quad\mathrm{with}\quad S_{0}=0.
\end{displaymath}
The probabilistic Lah numbers associated with $Y$ are given by
\begin{equation}
\frac{1}{k!}\bigg(E\bigg[\bigg(\frac{1}{1-t}\bigg)^{Y}\bigg]-1\bigg)^{k}=\sum_{n=k}^{\infty}L_{Y}(n,k)\frac{t^{n}}{n!},\quad (k\ge 0),\quad (\mathrm{see}\ [2,11,20,21]).\label{7}	
\end{equation}
When $Y=1$, $L_{Y}(n,k)=L(n,k),\ (n\ge k\ge 0)$. \par
From \eqref{7}, we note that
\begin{displaymath}
L_{Y}(n,k)=\frac{1}{k!}\sum_{l=0}^{k}\binom{k}{l}(-1)^{k-l}E\big[\langle S_{l}\rangle_{n}\big],\quad (\mathrm{see}\ [27]).
\end{displaymath}
In [24], the unsigned multi-Stirling numbers of the first kind are defined by
\begin{equation}
\mathrm{Li}_{k_{1},k_{2},\dots,k_{r}}(t)=\sum_{n=r}^{\infty}{n \brack k_{1},k_{2},\dots,k_{r}}\frac{t^{n}}{n!}.\label{8}
\end{equation}
From \eqref{8}, we note that
\begin{equation}
\sum_{n=r}^{\infty}{n \brack \underbrace{1,1,\dots,1}_{r-\mathrm{times}}}\frac{t^{n}}{n!}= \mathrm{Li}_{\underbrace{1,1,\dots,1}_{r-\mathrm{times}}}(t)=\frac{1}{r!}\big(-\log(1-t)\big)^{r}=\sum_{n=r}^{\infty}{n\brack r}\frac{t^{n}}{n!}.\label{9}	
\end{equation}
Thus, by \eqref{9}, we get
\begin{displaymath}
	{n \brack \underbrace{1,1,\dots,1}_{r-\mathrm{times}}}={n \brack r},\quad (n\ge r\ge 0).
\end{displaymath}
As the inverse relation of \eqref{9}, the multi-Stirling numbers of the second kind are given by
\begin{equation}
\mathrm{Li}_{k_{1},k_{2},\dots,k_{r}}\Big(1-e^{1-e^{t}}\Big)=\sum_{n=r}^{\infty}{n \brace k_{1},k_{2},\dots,k_{r}}\frac{t^{n}}{n!},\quad (\mathrm{see}\ [24]). \label{10}
\end{equation}
From \eqref{10}, we note that
\begin{equation}
\sum_{n=r}^{\infty}{n \brace \underbrace{1,1,\dots,1}_{r-\mathrm{times}}}\frac{t^{n}}{n!}=\mathrm{Li}_{\underbrace{1,1,\dots,1}_{r-\mathrm{times}}}\Big(1-e^{1-e^{t}}\Big)=\frac{1}{r!}\big(e^{t}-1\big)^{r}=\sum_{n=r}^{\infty}{n \brace r}\frac{t^{n}}{n!}.\label{11}
\end{equation}
Thus, by \eqref{11}, we get
\begin{equation}
{n \brace \underbrace{1,1,\dots,1}_{r-\mathrm{times}}}={n \brace r},\quad (n\ge r\ge 1). \label{12}
\end{equation}
By \eqref{9} and \eqref{10}, we get for $n\ge 0$ and $r\in\mathbb{N}$, we have
\begin{equation}
{n+1 \brace k_{1},k_{2},\dots,k_{r-1},1}=\sum_{m=r-1}^{n}\binom{n}{m}{m \brace k_{1},k_{2},\dots,k_{r-1}},\quad (\mathrm{see}\ [24]). \label{13}
\end{equation}
The probabilistic Stirling numbers of the second kind associated with $Y$ are defined by
\begin{equation}
{n \brace k}_{Y}=\frac{1}{k!}\sum_{j=0}^{k}\binom{k}{j}(-1)^{k-j}E\big[S_{j}^{n}\big],\quad (n\ge k\ge 0),\quad (\mathrm{see}\ [20-22,25,32]).\label{14}
\end{equation}
When $Y=1$,
\begin{displaymath}
	{n\brace k}_{Y}=\frac{1}{k!}\sum_{j=0}^{k}\binom{k}{j}(-1)^{k-j}j^{n}={n\brace k}.
\end{displaymath}

\section{Probabilistic multi-Stirling numbers of the second kind and probabilistic multi-Lah numbers }
For $k_{i}\ge 1$, $(i=1,2,\dots,r)$, by \eqref{1}, we get
\begin{align}
\frac{d}{dt} \mathrm{Li}_{k_{1},k_{2},\dots,k_{r}}(t)&=\frac{d}{dt}\sum_{0<m_{1}<m_{2}<\cdots<m_{r}}\frac{t^{m_{r}}}{m_{1}^{k_{1}}m_{2}^{k_{2}}\cdots m_{r}^{k_{r}}} \label{15} \\
&=\sum_{0<m_{1}<\cdots<m_{r}}\frac{t^{m_{r}-1}}{m_{1}^{k_{1}}m_{2}^{k_{2}}\cdots m_{r}^{k_{r}-1}}=\frac{1}{t} \mathrm{Li}_{k_{1},k_{2},\dots,k_{r-1},k_{r}-1}(t).\nonumber
\end{align}
Letting $k_{r}=1$ in \eqref{15}, we have
\begin{align}
\frac{d}{dt} \mathrm{Li}_{k_{1},k_{2},\dots,k_{r-1},1}(t)&=\sum_{0<m_{1}<m_{2}<\cdots<m_{r-1}}\frac{1}{m_{1}^{k_{1}}m_{2}^{k_{2}}\cdots m_{r-1}^{k_{r-1}}}\sum_{m_{r}=m_{r-1}+1}^{\infty}t^{m_{r}-1} \label{16} \\
&=\frac{1}{1-t}\sum_{0<m_{1}<\cdots<m_{r-1}}\frac{t^{m_{r-1}}}{m_{1}^{k_{1}}m_{2}^{k_{2}}\cdots m_{r-1}^{k_{r-1}}}=\frac{1}{1-t} \mathrm{Li}_{k_{1},k_{2},\dots,k_{r-1}}(t).\nonumber	
\end{align}
Thus, by \eqref{16}, we get
\begin{displaymath}
	\mathrm{Li}_{1,1}(t)=\frac{1}{2!}\big(-\log(1-t)\big)^{2},\quad \mathrm{Li}_{1,1,1}(t)=\frac{1}{3!}\big(-\log(1-t)\big)^{3}.
\end{displaymath}
Continuing this process, we have
\begin{equation}
\mathrm{Li}_{\underbrace{1,1,\dots,1}_{r-\mathrm{times}}}(t)=\frac{1}{r!}\big(-\log(1-t)\big)^{r},\quad (r\ge 1).\label{17}	
\end{equation}
Now, we define the {\it{probabilistic multi-Stirling numbers of the second kind associated with $Y$}} as
\begin{equation}
	\mathrm{Li}_{k_{1},k_{2},\dots,k_{r}}\Big(1-e^{1-E[e^{Yt}]}\Big)=\sum_{n=r}^{\infty}{n \brace k_{1},k_{2},\dots,k_{r}}_{Y}\frac{t^{n}}{n!}.\label{18}
\end{equation}
When $Y=1$, we have
\begin{displaymath}
{n \brace k_{1},k_{2},\dots,k_{r}}_{Y}={n \brace k_{1},k_{2},\dots,k_{r}}.
\end{displaymath}
Note that
\begin{align}
\sum_{n=r}^{\infty}{n \brace \underbrace{1,1,\dots,1}_{r-\mathrm{times}}}\frac{}{}_{Y}\frac{t^{n}}{n!}&=\mathrm{Li}_{\underbrace{1,1,\dots,1}_{r-\mathrm{times}}}\Big(1-e^{1-E[e^{Yt}]}\Big) \label{19}	\\
&=\frac{1}{r!}\Big(E\big[e^{Yt}\big]-1\Big)^{r}=\frac{1}{r!}\sum_{j=0}^{r}\binom{r}{j}(-1)^{r-j}E\big[e^{S_{j}t}\big] \nonumber\\
&=\sum_{n=0}^{\infty}\frac{1}{r!}\sum_{j=0}^{r}\binom{r}{j}(-1)^{r-j}E\big[S_{j}^{n}\big]\frac{t^{n}}{n!}.\nonumber
\end{align}
Thus, by \eqref{14} and \eqref{19}, we get
\begin{equation}
{n \brace \underbrace{1,1,\dots,1}_{r-\mathrm{times}}}\frac{}{}_{Y}={n \brace r}_{Y},\quad (n\ge r \ge 1).\label{20}
\end{equation}
From \eqref{1}, we note that
\begin{align}
\frac{d}{dt} \mathrm{Li}_{k_{1},k_{2},\dots,k_{r-1},1}&\Big(1-e^{1-E[e^{Yt}]}\Big)=\frac{d}{dt}\sum_{0<m_{1}<m_{2}<\cdots<m_{r}}\frac{(1-e^{1-E[e^{Yt}]})^{m_{r}}}{m_{1}^{k_{1}}m_{2}^{k_{2}}\cdots m_{r-1}^{k_{r-1}}m_{r}} \label{21} \\
&=\sum_{0<m_{1}<\cdots<m_{r-1}}\frac{E[Ye^{Yt}]e^{1-E[e^{Yt}]}}{m_{1}^{k_{1}}\cdots m_{r-1}^{k_{r-1}}}\sum_{n_{r}=m_{r-1}+1}^{\infty}\Big(1-e^{1-E[e^{Yt}]}\Big)^{m_{r}-1} \nonumber \\
&=\sum_{0<m_{1}<\cdots<m_{r-1}}\frac{E[Ye^{Yt}]e^{1-E[e^{Yt}]}}{m_{1}^{k_{1}}\cdots m_{r-1}^{k_{r-1}}} \frac{(1-e^{1-E[e^{Yt}]})^{m_{r-1}}} {e^{1-E[e^{Yt}]}}\nonumber  \\
&=E\big[Ye^{Yt}\big] \mathrm{Li}_{k_{1},k_{2},\dots,k_{r-1}}\Big(1-e^{1-E[e^{Yt}]}\Big)\nonumber \\
&=\sum_{l=0}^{\infty}E\big[Y^{l+1}\big]\frac{t^{l}}{l!}\sum_{m=r-1}^{\infty}{m \brace k_{1},\dots,k_{r-1}}_{Y}\frac{t^{m}}{m!}\nonumber \\
&=\sum_{n=r-1}^{\infty}\sum_{m=r-1}^{n}\binom{n}{m}{m \brace k_{1},k_{2},\dots,k_{r-1}}_{Y}E\big[Y^{n-m+1}\big]\frac{t^{n}}{n!}.\nonumber
\end{align}
On the other hand, by \eqref{18}, we get
\begin{align}
\frac{d}{dt} \mathrm{Li}_{k_{1},k_{2},\dots,k_{r-1},1}\Big(1-e^{1-E[e^{Yt}]}\Big)&=\frac{d}{dt}\sum_{n=r}^{\infty}{n \brace k_{1},k_{2},\dots,k_{r-1},1}_{Y}\frac{t^{n}}{n!} \label{22} \\
&=\sum_{n=r}^{\infty}{n \brace k_{1},k_{2},\dots,k_{r-1},1}_{Y}\frac{t^{n-1}}{(n-1)!}\nonumber\\
&=\sum_{n=r-1}^{\infty}{n+1 \brace k_{1},k_{2},\dots,k_{r-1},1}_{Y}\frac{t^{n}}{n!}. \nonumber
\end{align}
Therefore, by \eqref{20}, \eqref{21} and \eqref{22}, we obtain the following theorem,
\begin{theorem}
For $n,r\in\mathbb{N}$ with $n\ge r-1$, we have
\begin{displaymath}
\sum_{m=r-1}^{n}\binom{n}{m}E\big[Y^{n-m+1}\big]{m \brace k_{1},k_{2},\dots,k_{r-1}}_{Y}={n+1 \brace k_{1},k_{2},\dots,k_{r-1},1}_{Y}.
\end{displaymath}
In particular, for $n \ge r$, we have
\begin{equation*}
\sum_{m=r-1}^{n-1}\binom{n-1}{m}{m \brace r-1}_{Y}E[Y^{n-m}]={n\brace r}_{Y},
\quad \sum_{m=r-1}^{n-1}\binom{n-1}{m}{m \brace r-1}={n\brace r}.
\end{equation*}
\end{theorem}
Replacing $t$ by $E[e^{Yt}]-1$ in \eqref{4}, on the one hand, we get
\begin{align}
\frac{1}{(1-e^{1-E[e^{Yt}]})^{r}} \mathrm{Li}_{k_{1},k_{2},\dots,k_{r-1},k_{r}}\Big(1-e^{1-E[e^{Yt}]}\Big)&=\sum_{m=0}^{\infty}B_{m}^{(k_{1},\dots,k_{r})}\frac{1}{m!}\Big(E[e^{Yt}]-1\Big)^{m} \label{23} \\
&=\sum_{m=0}^{\infty}B_{m}^{(k_{1},k_{2},\dots,k_{r})}\sum_{n=m}^{\infty}{n \brace m}_{Y}\frac{t^{n}}{n!}\nonumber\\
&=\sum_{n=0}^{\infty}\bigg(\sum_{m=0}^{n}B_{m}^{(k_{1},k_{2},\dots,k_{r})}{n\brace m}_{Y}\bigg)\frac{t^{n}}{n!}. \nonumber
\end{align}
On the other hand, by \eqref{18}, we get
\begin{align}
&\frac{1}{(1-e^{1-E[e^{Yt}]})^{r}} \mathrm{Li}_{k_{1},k_{2},\dots,k_{r-1},k_{r}}\Big(1-e^{1-E[e^{Yt}]}\Big)\label{24}\\
&=\sum_{l=0}^{\infty}\binom{r+l-1}{l}e^{l(1-E[e^{Yt}])}\sum_{m=r}^{\infty}{m \brace k_{1},k_{2},\dots,k_{r}}_{Y}\frac{t^{m}}{m!} \nonumber	\\
&=\sum_{l=0}^{\infty}\binom{r+l-1}{l}\sum_{j=0}^{\infty}l^{j}(-1)^{j}\frac{1}{j!}\Big(E[e^{Yt}]-1\Big)^{j}\sum_{m=r}^{\infty}{m \brace k_{1},k_{2},\dots,k_{r}}_{Y}\frac{t^{m}}{m!}\nonumber \\
&=\sum_{l=0}^{\infty}\binom{r+l-1}{l}\sum_{j=0}^{\infty}l^{j}(-1)^{j}\sum_{i=j}^{\infty}{i \brace j}_{Y}\frac{t^{i}}{i!}\sum_{m=r}^{\infty}{m \brace k_{1},k_{2},\dots,k_{r}}\frac{t^{m}}{m!}\nonumber \\
&=\sum_{l=0}^{\infty}\binom{r+l-1}{l}\sum_{i=0}^{\infty}\bigg(\sum_{j=0}^{i}l^{j}(-1)^{j}{i \brace j}_{Y}\frac{t^{i}}{i!}\bigg)\sum_{m=r}^{\infty}{m \brace k_{1},k_{2},\dots,k_{r}}_{Y}\frac{t^{m}}{m!}\nonumber \\
&=\sum_{n=r}^{\infty}\bigg(\sum_{l=0}^{\infty}\sum_{i=0}^{n-r}\sum_{j=0}^{i}\binom{n}{i}\binom{r+l-1}{l}l^{j}(-1)^{j}{i \brace j}_{Y}{n-i \brace k_{1},k_{2},\dots,k_{r}}_{Y}\bigg)\frac{t^{n}}{n!}.\nonumber
\end{align}
Therefore, by \eqref{23} and \eqref{24}, we obtain the following theorem.
\begin{theorem}
For $n,r\in\mathbb{N}$ with $n\ge r$, we have
\begin{equation*}
\sum_{m=0}^{n}B_{m}^{(k_{1},k_{2},\dots,k_{r})}{n \brace m}_{Y}
=\sum_{l=0}^{\infty}\sum_{i=0}^{n-r}\sum_{j=0}^{i}\binom{n}{i}\binom{r+l-1}{l}l^{j}(-1)^{j}{i \brace j}_{Y}{n-i \brace k_{1},k_{2},\dots,k_{r}}_{Y}.
\end{equation*}
\end{theorem}
Replacing $t$ by $1-e^{1-E[e^{Yt}]}$ in \eqref{8}, we get
\begin{align}
\mathrm{Li}_{k_{1},k_{2},\dots,k_{r}}\Big(1-e^{1-E[e^{Yt}]}\Big)&=\sum_{m=r}^{\infty}{m \brack k_{1},k_{2},\dots,k_{r}}\frac{1}{m!}\Big(1-e^{1-E[e^{Yt}]}\Big)^{m} \label{25} \\
&=\sum_{m=r}^{\infty}{m \brack k_{1},k_{2},\dots,k_{r}}(-1)^{m}\sum_{l=m}^{\infty}{l \brace m}\frac{1}{l!}\Big(1-E\big[e^{Yt}\big]\Big)^{l} \nonumber \\
&=\sum_{l=r}^{\infty}\sum_{m=r}^{l}{m \brack k_{1},k_{2},\dots,k_{r}}(-1)^{m-l}{l \brace m}\sum_{n=l}^{\infty}{n\brace l}_{Y}\frac{t^{n}}{n!}\nonumber \\
&=\sum_{n=r}^{\infty}\sum_{m=r}^{l}{m \brack k_{1},k_{2},\dots,k_{r}}(-1)^{m-l}{l \brace m}\sum_{l=r}^{n}{n \brace l}_{Y}\frac{t^{n}}{n!}\nonumber \\
&=\sum_{n=r}^{\infty}\bigg(\sum_{l=r}^{n}\sum_{m=r}^{l}{m \brack k_{1},k_{2},\dots,k_{r}}(-1)^{m-l}{l \brace m}{n \brace l}_{Y}\bigg)\frac{t^{n}}{n!}.\nonumber
\end{align}
Therefore ,by \eqref{18} and \eqref{25}, we obtain the following theorem.
\begin{theorem}
For $n,r\in\mathbb{N}$ with $n\ge r$, we have
\begin{displaymath}
{n \brace k_{1},k_{2},\dots,k_{r}}_{Y}=\sum_{l=r}^{n}\sum_{m=r}^{l}(-1)^{m-l}{l \brace m}{n \brace l}_{Y}{m \brack k_{1},k_{2},\dots,k_{r}}.
\end{displaymath}
\end{theorem}
Now, we consider the {\it{probabilistic multi-Lah numbers associated with $Y$}} which are given by
\begin{equation}
\mathrm{Li}_{k_{1},k_{2},\dots,k_{r}}\bigg(1-e^{1-E\big[\big(\frac{1}{1-t}\big)^{Y}\big)\big]}\bigg)=\sum_{n=r}^{\infty}L_{Y}^{(k_{1},k_{2},\dots,k_{r})}(n,r)\frac{t^{n}}{n!}.\label{26}
\end{equation}
By \eqref{17}, we note that
\begin{align}
\sum_{n=r}^{\infty}L_{Y}^{\overbrace{(1,1,\dots,1)}^{r-\mathrm{times}}}(n,r)\frac{t^{n}}{n!}&=\mathrm{Li}_{\underbrace{1,1,\dots,1}_{r-\mathrm{times}}}\bigg(1-e^{1-E\big[\big(\frac{1}{1-t}\big)^{Y}\big]}\bigg)\label{27} \\
&=\frac{1}{r!}\bigg(E\bigg[\bigg(\frac{1}{1-t}\bigg)^{Y}\bigg]-1\bigg)^{r}\nonumber\\
&=\sum_{n=r}^{\infty}L_{Y}(n,r)\frac{t^{n}}{n!}. \nonumber
\end{align}
Thus, by \eqref{27}, we get
\begin{equation}
L_{Y}^{\overbrace{(1,1,\dots,1)}^{r-\mathrm{times}}}(n,r)=L_{Y}(n,r),\quad (n\ge r>0).\label{28}
\end{equation}
Using \eqref{26} and \eqref{18} with $t$ replaced by$-\log(1-t)$, we get
\begin{align}
\sum_{n=r}^{\infty} L_{Y}^{(k_{1},k_{2},\dots,k_{r})}(n,r)\frac{t^{n}}{n!}&= \mathrm{Li}_{k_{1},k_{2},\dots,k_{r}}\Big(1-e^{1-E\big[\big(\frac{1}{1-t}\big)^{Y}\big]}\Big)\label{29} \\
&=\sum_{k=r}^{\infty}{k \brace k_{1},k_{2},\dots,k_{r}}_{Y}\frac{1}{k!}\big(-\log(1-t)\big)^{k}\nonumber \\
&=\sum_{k=r}^{\infty}{n \brace k_{1},k_{2},\dots,k_{r}}_{Y}\sum_{n=k}^{\infty}{n \brack k}\frac{t^{n}}{n!}\nonumber \\
&=\sum_{n=r}^{\infty}\sum_{k=r}^{n}{n \brace k_{1},k_{2},\dots,k_{r}}_{Y}{n \brack k}\frac{t^{n}}{n!}. \nonumber
\end{align}
Therefore, by comparing the coefficients on both sides of \eqref{29}, we obtain the following theorem.
\begin{theorem}
For $n,r\in\mathbb{N}$ with $n\ge r$, we have
\begin{displaymath}
L_{Y}^{(k_{1},k_{2},\dots,k_{r})}(n,r)= \sum_{k=r}^{n}{n \brace k_{1},k_{2},\dots,k_{r}}_{Y}{n \brack k}.
\end{displaymath}
\end{theorem}
Now, from \eqref{4} and \eqref{18}, we observe that
\begin{align}
\sum_{n=r}^{\infty}&{n \brace k_{1},k_{2},\dots,k_{r}}_{Y}\frac{t^{n}}{n!}= \mathrm{Li}_{k_{1},k_{2},\dots,k_{r}}\Big(1-e^{1-E[e^{Yt}]}\Big)\label{30} \\
&=\frac{\mathrm{Li}_{k_{1},k_{2},\dots,k_{r}}\big(1-e^{1-E[e^{Yt}]}\big)}{\big(1-e^{1-E[e^{Yt}]}\big)^{r}}\Big(1-e^{1-E[e^{Yt}]}\Big)^{r}\nonumber \\
&=\sum_{m=0}^{\infty}B_{m}^{(k_{1},k_{2},\dots,k_{r})}\frac{1}{m!}\Big(E[e^{Yt}]-1\Big)^{m}(-1)^{r}r!\frac{1}{r!}\Big(e^{1-E[e^{Yt}]}-1\Big)^{r} \nonumber\\
&=\sum_{m=0}^{\infty}B_{m}^{(k_{1},k_{2},\dots,k_{r})}\frac{1}{m!}\Big(E[e^{Yt}]-1\Big)^{m}(-1)^{r}r!\sum_{l=r}^{\infty}{l \brace r}\frac{1}{l!}\Big(1-E[e^{Yt}]\Big)^{l} \nonumber \\
&=\sum_{m=0}^{\infty}\sum_{l=r}^{\infty}B_{m}^{(k_{1},k_{2},\dots,k_{r})}r!(-1)^{l-r}{l \brace r}\frac{(m+l)!}{m!l!}\frac{1}{(m+l)!}\Big(E\big[e^{Yt}\big]-1\Big)^{m+l} \nonumber \\
&=\sum_{m=0}^{\infty}\sum_{l=r}^{\infty}B_{m}^{(k_{1},k_{2},\dots,k_{r})}r!(-1)^{l-r}{l \brace r}\binom{m+l}{m}\sum_{n=m+l}^{\infty}{n \brace m+l}_{Y}\frac{t^{n}}{n!}\nonumber \\
&=\sum_{n=r}^{\infty}\bigg(\sum_{m=0}^{n-r}\sum_{l=r}^{n-m}B_{m}^{(k_{1},k_{2},\dots,k_{r})}r!(-1)^{l-r}\binom{m+l}{m}{l \brace r}{n \brace m+l}_{Y}\bigg)\frac{t^{n}}{n!}. \nonumber
\end{align}
By comparing the coefficients on both sides of \eqref{30}, we obtain the following theorem.
\begin{theorem}
Let $n,r$ be positive integers with $n\ge r$. Then we have
\begin{displaymath}
{n \brace k_{1},k_{2},\dots,k_{r}}_{Y}=\sum_{m=0}^{n-r}\sum_{l=r}^{n-m}r!(-1)^{l-r}\binom{m+l}{m}{l \brace r}{n \brace m+l}_{Y}B_{m}^{(k_{1},k_{2},\dots,k_{r})}.
\end{displaymath}
In particular, we have
\begin{displaymath}
{n \brace r}_{Y}=\sum_{m=0}^{n-r}\sum_{l=r}^{n-m}(-1)^{m+l-r}\binom{m+l}{m}{l \brace r}{n \brace m+l}_{Y}B_{m}^{(r)}.
\end{displaymath}
\end{theorem}
For $r\in\mathbb{N}$, the probabilistic Fubini polynomials of order $r$ associated with $Y$ are defined by
\begin{equation}
\frac{1}{\big(1-y(E[e^{Yt}]-1)\big)^{r}}=\sum_{n=0}^{\infty}F_{n}^{(r,Y)}(y)\frac{t^{n}}{n!},\quad(\mathrm{see}\ [32]).\label{31}
\end{equation}
Replacing $t$ by $E[e^{Yt}]-1$ in \eqref{2}, on the one hand, we have
\begin{align}
\frac{1}{(2-E[e^{Yt}])^{r}}\mathrm{Li}_{k_{1},k_{2},\dots,k_{r}}\Big(1-e^{1-E[e^{Yt}]}\Big)&=\sum_{k=r}^{\infty}L^{(k_{1},k_{2},\dots,k_{r})}(k,r)\frac{1}{k!}\Big(E\big[e^{Yt}\big]-1\Big)^{k}\label{32} \\
&=\sum_{k=r}^{\infty}L^{(k_{1},k_{2},\dots,k_{r})}(k,r)\sum_{n=k}^{\infty}{n \brace k}_{Y}\frac{t^{n}}{n!}\nonumber\\
&=\sum_{n=r}^{\infty}\bigg(\sum_{k=r}^{n}{n \brace k}_{Y}L^{(k_{1},k_{2},\dots,k_{r})}(k,r)\bigg)\frac{t^{n}}{n!}.\nonumber	
\end{align}
On the other hand, by \eqref{18} and \eqref{31}, we get
\begin{align}
\frac{1}{(2-E[e^{Yt}])^{r}}\mathrm{Li}_{k_{1},k_{2},\dots,k_{r}}\Big(1-e^{1-E[e^{Yt}]}\Big)&=\sum_{k=r}^{\infty}{k \brace k_{1},k_{2},\dots,k_{r}}_{Y}\frac{t^{k}}{k!}\sum_{l=0}^{\infty}F_{l}^{(r,Y)}(1)\frac{t^{l}}{l!}\label{33} \\
&=\sum_{n=r}^{\infty}\bigg(\sum_{k=r}^{n}{k \brace k_{1},k_{2},\dots,k_{r}}_{Y}\binom{n}{k}F_{n-k}^{(r,Y)}(1)\bigg)\frac{t^{n}}{n!}.\nonumber
\end{align}
Therefore, by \eqref{32} and \eqref{33}, we obtain the following theorem.
\begin{theorem}
For $n,r\in\mathbb{N}$ with $n\ge r$, we have
\begin{displaymath}
\sum_{k=r}^{n}{n \brace k}_{Y}L^{(k_{1},k_{2},\dots,k_{r})}(k,r)=\sum_{k=r}^{n}\binom{n}{k}{k \brace k_{1},k_{2},\dots,k_{r}}_{Y}F_{n-k}^{(r,Y)}(1).
\end{displaymath}
\end{theorem}

\section{Conclusion}
Let $Y$ be a random variable whose moment generating function exists in a neighborhood of the origin. The Stirling numbers of the second ${n \brace r}$ are generalized to the multi-Stirling numbers of the second kind ${n \brace k_{1},k_{2},\dots,k_{r}}$ and further to the probabilistic multi-Stirling numbers of the second kind associated with $Y$, ${n \brace k_{1},k_{2},\dots,k_{r}}_{Y}$. Likewise, the unsigned Lah numbers $L(n,r)$ are generalized to the multi-Lah numbers $L^{(k_{1},k_{2},\dots,k_{r})}(n,r)$ and further to the probabilistic multi-Lah numbers associated with $L_{Y}^{(k_{1},k_{2},\dots,k_{r})}(n,r)$. In connection with several other special numbers, we investigated some properties, related identities, recurrence relations and explicit expressions of the numbers ${n \brace k_{1},k_{2},\dots,k_{r}}_{Y}$ and $L_{Y}^{(k_{1},k_{2},\dots,k_{r})}(n,r)$. \par
In recent years, degenerate versions, $\lambda$-analogues and probabilistic extensions of many special numbers and polynomials have been studied and many remarkable results have been obtained (see [2,3,10,11,16-18,20,21,23,25-28,32] and the references therein). It is one of our future research projects to continue to study degenerate versions, $\lambda$-analogues and probabilistic extensions of various special numbers and polynomials and to find their applications to physics, science and engineering as well as to mathematics.

\vspace{0.1in}

{\bf{Declarations}}

\vspace{0.1in}

{\bf{Ethics approval and consent to participate}} The submitted article is an original research paper and has not been published anywhere else.

\vspace{0.1in}

{\bf{Consent for publication}} All authors consent for publication.

\vspace{0.1in}

{\bf{Data Availability Statement}} There is no data used in this article.

\vspace{0.1in}

{\bf{Competing interests}} The authors declare no competing interests.
\vspace{0.1in}

{\bf{Funding}}
This research has been conducted by the Research Grant of Kwangwoon University in 2024.

\vspace{0.1in}

{\bf{Authors' contributions}} All authors contributed equally in the preparation of this manuscript.


\begin{thebibliography}{9}
\bibitem{1}
Abramowitz, M. and Stegun, I. A. \emph{Handbook of Mathematical Functions with Formulas, Graphs, and Mathematical Tables,} National Bureau of Standards Applied Mathematics Series, U. S. Government Printing Office, Washington, DC, 1964.
\bibitem{2}
Adell, J. A.; B\'enyi, B. \emph{Probabilistic stirling numbers and applications,} Aequat. Math. (2024). https://doi.org/10.1007/s00010-024-01073-1
\bibitem{3}
Aydin, M. S.; Acikgoz, M.; Araci, S. \emph{A new construction on the degenerate Hurwitz-zeta function associated with certain applications,} Proc. Jangjeon Math. Soc. \textbf{25} (2022), no. 2, 195-203.
\bibitem{4}
Bening, V. E. \emph{On the asymptotic deficiency of some statistical estimators based on samples with random size,} Proc. Jangjeon Math. Soc. \textbf{21} (2018), no. 2, 185-193.
\bibitem{5}
Blyth, C. R.; Pathak, P. K. Notes: \emph{A note on easy proofs of Stirling's theorem,} Amer. Math. Monthly \textbf{93} (1986), no. 5, 376-379.
\bibitem{6}
Boubellouta, K.; Boussayoud, A.; Araci, S.; Kerada, M. \emph{Some theorems on generating functions and their applications,} Adv. Stud. Contemp. Math. (Kyungshang) \textbf{30}(2020), no. 3, 307-324.
\bibitem{7}
Boyadzhiev, K. N. \emph{Convolutions for Stirling numbers, Lah numbers, and binomial coefficients,} Proc. Jangjeon Math. Soc. \textbf{25} (2022), no. 2, 227-244.
\bibitem{8}
Broder, A. Z. \emph{The $r$-Stirling numbers,} Discrete Math. \textbf{49} (1984), 241-259.
\bibitem{9}
Carlitz, L. \emph{Degenerate Stirling, Bernoulli and Eulerian numbers,} Utilitas Math. \textbf{15} (1979), 51-88.
\bibitem{10}
Chen, L.; Dolgy, D. V.; Kim, T.; Kim, D. S. \emph{Probabilistic type 2 Bernoulli and Euler polynomials,} AIMS Math. \textbf{9} (2024), no. 6, 14312-14324.
\bibitem{11}
Chen, L.; Kim, T.; Kim, D. S.; Lee, H.; Lee, S.-H. \emph{Probabilistic degenerate central Bell polynomials,} Math. Comput. Model. Dyn. Syst., https://doi.org/10.1080/13873954.2024.2358899.
\bibitem{12}
Comtet, L. \emph{Advanced Combinatorics. The Art of Finite and Infinite Expansions,} Revised and enlarged edition., D. Reidel Publishing Co., Dordrecht, 1974.
\bibitem{13}
Graham, R. L.; Knuth, D. E.; Patashnik, O. \emph{Concrete Mathematics: A Foundation for Computer Science,} 2 Eds., Addison Wesley Publishing Company, Massachusetts, 1994.
\bibitem{14}
Gun, D.; Simsek, Y. \emph{Combinatorial sums involving Stirling, Fubini, Bernoulli numbers and approximate values of Catalan numbers,} Adv. Stud. Contemp. Math. (Kyungshang) \textbf{30} (2020), no. 4, 503-513.
\bibitem{15}
Hwang, K.-S. \emph{Almost sure convergence of weighted sums for widely negative dependent random variables under sub-linear expectations,} Adv. Stud. Contemp. Math. (Kyungshang) \textbf{32} (2022), no. 4, 465-477.
\bibitem{16}
Jang, L.-C. \emph{A note on degenerate type 2 multi-poly-Genocchi polynomials,} Adv. Stud. Contemp. Math. (Kyungshang) \textbf{30} (2020), no. 4, 537-543.
\bibitem{17}
Kim, D. S.; Kim, T. \emph{Normal ordering associated with $\lambda$-Whitney numbers of the first Kind in $\lambda$-shift algebra,}  Russ. J. Math. Phys. \textbf{30} (2023), no. 3, 310-319.
\bibitem{18}
Kim, D. S.; Kim, T. \emph{A note on a new type of degenerate Bernoulli numbers,} Russ. J. Math. Phys. \textbf{27} (2020), no. 2, 227-235.
\bibitem{19}
Kim, M.-S.; Kim, T. \emph{An explicit formula on the generalized Bernoulli number with order $n$,} Indian J. Pure Appl. Math. \textbf{31} (2000), no. 11, 1455-1461.
\bibitem{20}
Kim, T.; Kim, D. S. \emph{Probabilistic Bernoulli and Euler polynomials,} Russ. J. Math. Phys. \textbf{31} (2024), no. 1, 94-105.
\bibitem{21}
Kim, T.; Kim, D. S. \emph{Probabilistic degenerate Bell polynomials associated with random variables,} Russ. J. Math. Phys. \textbf{30} (2023), no. 4, 528-542.
\bibitem{22}
Kim, T.; Kim, D. S. \emph{Generalization of Spivey's recurrence relation},  Russ. J. Math. Phys. \textbf{31} (2023), no. 2, 218-226.
\bibitem{23}
Kim, T.; Kim, D. S. \emph{Some identities involving degenerate Stirling numbers associated with several degenerate polynomials and numbers,} Russ. J. Math. Phys. \textbf{30} (2023), no. 1, 62-75.
\bibitem{24}
Kim, T.; Kim, D. S.; Kim, H. K. \emph{Multi-Stirling numbers of the second kind} arXiv:2303.00275
\bibitem{25}
Kim, T.; Kim, D. S.; Kwon, J. \emph{Probabilistic degenerate Stirling polynomials of the second kind and their applications.} Math. Comput. Model. Dyn. Syst. \textbf{30} (2024), no. 1, 16-30.
\bibitem{26}
Luo, L.; Kim, T.; Kim, D. S.; Ma, Y. \emph{Probabilistic degenerate Bernoulli and degenerate Euler polynomials,} Math. Comput. Model. Dyn. Syst. \textbf{30} (2024), no. 1, 342-363.
\bibitem{27}
Ma. Y.; Kim, T.; Kim, D. S. \emph{Probabilistic Lah numbers and Lah-Bell polynomials,} arXiv:2406.01200
\bibitem{28}
Park, J.-W. \emph{On the degenerate multi-poly-Genocchi polynomials and numbers,} Adv. Stud. Contemp. Math. (Kyungshang) \textbf{33} (2023), no. 2, 181-186.
\bibitem{29}
Roman, S. \emph{The umbral calculus,} Pure and Applied Mathematics, 111, Academic Press, Inc. [Harcourt Brace Jovanovich, Publishers], New York, 1984.
\bibitem{30}
Ross, S. M. \emph{Introduction to probability models,} Thirteenth edition, Academic Press, London, 2024.
\bibitem{31}
Simsek, Y. \emph{Identities and relations related to combinatorial numbers and polynomials,} Proc. Jangjeon Math. Soc. \textbf{20} (2017), no. 1, 127-135.
\bibitem{32}
Soni, R.; Pathak, A.K.; Vellaisamy, P. \emph{A probabilistic extension of the Fubini polynomials,} Bull. Malays. Math. Sci. Soc. \textbf{47}, 102 (2024).
\end{thebibliography}
\end{document}